\documentclass{elsart}

\usepackage[T1]{fontenc}
\usepackage{amsmath,amsfonts,amssymb}
\usepackage{mathrsfs} 
\usepackage{enumerate}

\makeatletter
\newenvironment{jbRemark}{\stepcounter{thm}\let\@currentlabel\thethm\emph{Remark \thethm.}}
\makeatother

\newenvironment{jbProof}{\textbf{Proof.}\ }{\qed}

\newcounter{jbCmptAssA}
\newcounter{jbCmptAssB}
\renewcommand{\thejbCmptAssA}{$\mathbf{A_{\arabic{jbCmptAssA}}}$}
\newenvironment{assumptions}
  {\begin{list}{(\thejbCmptAssA)}{\setlength{\listparindent}{0pt}\usecounter{jbCmptAssA}\setcounter{jbCmptAssA}{\value{jbCmptAssB}}}}
  {\end{list}\setcounter{jbCmptAssB}{\value{jbCmptAssA}}}

\newcommand{\jbNoNum}{}

\DeclareMathOperator{\diverg}{div}                                
\DeclareMathOperator{\jacob}{Jac}
\DeclareMathOperator{\supp}{supp}

\newcommand{\scal}[2]{{\left\langle \,{#1}, \,{#2} \,\right\rangle}}
\newcommand{\scalfct}{{\left\langle\cdot,\cdot\right\rangle}}
\newcommand{\Pro}[1]{\mathbb{P}_{#1}}

\newcommand{\EspLetter}{\mathbb{E}}

\newcommand{\Espx}[2]{\EspLetter_{#1} \left\{ #2 \right\}}
\newcommand{\F}{\mathcal{F}}

\newcommand{\Cclass}{C}
\newcommand{\Ck}[2]{{\Cclass}^{#1}({#2})}
\newcommand{\Ckc}[2]{{\Cclass}_{c}^{#1}({#2})}
\newcommand{\dd}{k}
\newcommand{\SDE}{\mathscr{S}}
\newcommand{\Ito}{It\=o\ }
\newcommand{\dDiff}{\,\mathrm{d}}
\newcommand{\dStrato}{ \circ \dDiff }
\newcommand{\dX}{\dDiff X}
\newcommand{\dt}{\dDiff t}
\newcommand{\dBIto}{\dDiff B}
\newcommand{\dBStrato}{\dStrato B}
\newcommand{\dm}{\dDiff m}
\newcommand{\ds}{\dDiff s}
\newcommand{\dsigma}{\dDiff \sigma}
\newcommand{\dpdf}{\dDiff p}
\newcommand{\VF}{A}
\newcommand{\Jout}{{J}^{\mathrm{out}}}
\newcommand{\Jin}{{J}^{\mathrm{in}}}
\newcommand{\tauEtoile}{\tau^{*}}   
\newcommand{\indsaut}{j}
\newcommand{\OpKolBack}{L}        
\newcommand{\OpKolForw}{\OpKolBack^{*}}
\newcommand{\Loi}[1]{{\mu}_{#1}}
\newcommand{\M}{\mathbb{M}}
\newcommand{\Mb}{\overline{\M}}
\newcommand{\dM}{\partial\M}
\newcommand{\rvm}{m}
\newcommand{\St}{\mathbb{S}}
\newcommand{\Stb}{\overline{\St}}
\newcommand{\T}{\mathbb{T}}
\newcommand{\Hyp}{H}
\newcommand{\MHyp}{\M\setminus\Hyp}
\newcommand{\Gua}{G}

\begin{document}


\begin{frontmatter}
  
  \title{Fokker-Planck-Kolmogorov equation\\ for stochastic
    differential equations\\ with boundary hitting resets}
  
  \author{Julien Bect\corauthref{cor}\thanksref{email}}, 
  \author{Hana Baili\thanksref{email}},
  \author{Gilles Fleury\thanksref{email}}

  \corauth[cor]{Corresponding author.}
  \thanks[email]{Email adresses: Firstname.Lastname@Supelec.fr}

  \address{Department of Signal Processing and Electronic Systems,\\
   \'{E}cole Sup\'{e}rieure d'\'{E}lectricit{\'e},\\
   3 rue Joliot-Curie, Plateau de Moulon,\\
   91192 Gif-sur-Yvette cedex, France.}
  
  \begin{abstract} We consider a Markov process on a Riemannian
    manifold, which solves a stochastic differential equation in the
    interior of the manifold and jumps according to a deterministic
    reset map when it reaches the boundary.  We derive a partial
    differential equation for the probability density function,
    involving a non-local boundary condition which accounts for the
    jumping behaviour of the process. This is a generalisation of the
    usual Fokker-Planck-Kolmogorov equation for diffusion processes.
    The result is illustrated with an example in the field of
    stochastic hybrid systems.
  \end{abstract}
  
  \begin{keyword}
    Stochastic differential equations \sep Fokker-Planck equation \sep
    Kolmogorov's forward equation \sep Stochastic Hybrid Systems \MSC
    58J65 \sep 60H10 \sep 60J60 \sep 60J75
  \end{keyword}

\end{frontmatter}

\section{Introduction}

This paper investigates Kolmogorov's forward equation, also called the
Fokker-Planck equation (FPE), for a class of Markov processes whose
behaviour can be described as follows: let $\Mb$ be a $n$-dimensional
smooth\footnote{``smooth'' stands for $\Cclass^\infty$, here and throughout
  the whole paper, even though we do not really need that much
  regularity to prove our results.} manifold with boundary $\dM$ and
interior $\M$, and let $\St = \M \cup \T$ be the topological space
obtained by adding to $\M$ a finite (possibly empty) set $\T$ of
isolated points. Our process $X$ evolves in $\M$ according to a
Stratonovich stochastic differential equation\footnote{see e.g.
  Ikeda-Watanabe \cite[Chapters 3 and 5]{IW89} for the basic
  definitions.} (SDE) on $\M$,
\begin{equation}
  \tag{$\SDE$}\label{eq:SDE}
  \dX_t \;=\; \VF_0 (X_t) \dt \;+\; \sum_{r=1}^{n}\, \VF_r(X_t)\, \dBStrato_t^{r} \,,
\end{equation}
where $B$ is a $n$-dimensional Wiener process, $\dStrato$ denotes the
Stratonovich differential, and the $\VF_r$'s are $n+1$ smooth vector
fields on $\Mb$. If $X$ happens to hit the boundary $\dM$, it is
instantaneously reset either in $\M$ or to one of the terminal states
$x \in \T$, according to a measurable map $\Phi:\dM\to\St$, that we call the
reset map.

Our interest for such processes has been motivated by a class of
continuous-time Markovian models called Stochastic Hybrid Systems
(SHS), introduced recently by Hu et al. in \cite{Hu:2000:TTSHS}.
Roughly speaking, a SHS is obtained by replacing the ordinary
differential equations by SDEs in a deterministic hybrid system. SHSs
are not, of course, the only way to introduce randomness into the
deterministic framework of hybrid dynamical systems. Notably, Davis
\cite{Davis:1984:PDMP} has defined and studied a wide class of
non-diffusive hybrid processes. See \cite{Pola:2003:StoHybModels} and
the references therein for a survey of the litterature on
continuous-time stochastic hybrid models.

The link between diffusion processes and the FPE has been known since
the early days of Markov processes
\cite{Kolmogorov:1931:Analytischen}. Let us recall some basic facts
concerning the case of Stratonovich SDEs on a manifold without
boundary. From now on, $\Mb$ is endowed with a Riemannian metric
$\scalfct$ and we denote by $\rvm$ the Riemannian volume measure on
$\M$. Assume that $\dM=\varnothing$ and that the solution of
\eqref{eq:SDE} has a density $p_t$ with respect to $\rvm$, for all $t \geq
0$. Then, under some additional smoothness assumptions, $p_t$ solves
the partial differential equation (PDE)
\begin{equation}
  \frac{\partial p_t}{\partial t} \;=\; \OpKolForw p_t \,,
  \jbNoNum
\end{equation}
where 
\begin{equation}
  \label{eq:Hormander}
  \OpKolBack \;=\; \VF_0 \,+\, \frac{1}{2}\, \sum_{r=1}^n\, \VF_r^2 
\end{equation}
is the infinitesimal generator of the process, written in H{\"o}rmander
form\footnote{the notation $\VF_r^2$\ stands for the second order
  differential operator ${\VF_r}\circ{\VF_r}$, i.e.  in coordinates
  $\VF_r^2 \varphi = \sum_{i,j} \VF_r^i\ \partial_i ( \VF_r^j \ \partial_j \varphi )$.}, and
$\OpKolForw$ denotes its formal adjoint. Introducing the probability
current vector
\begin{equation}
  \label{eq:Jdef}
  J_t \;=\; p_t\, \VF_0 \;-\; \frac{1}{2}\, \sum_{r=1}^{n}\, \diverg \left( p_t \VF_r \right)\, \VF_r \,,
\end{equation}
where $\diverg$ denotes the divergence operator on $\M$,
the FPE can be rewritten as a local conservation
equation \cite{Strato:1966:StochIntegrals,Gardiner:1985:Handbook}:
\begin{equation}
  \label{eq:LocalConservation}
  \frac{\partial p_t}{\partial t} + \diverg (J_t) \;=\; 0 \,.
\end{equation}

The point of this paper is to prove that, for a SDE with boundary
hitting resets, under sufficient assumptions:
\begin{enumerate}[-]
\item Equation \eqref{eq:LocalConservation} holds in $\MHyp$, where
  $\Hyp = \Phi \left(\dM\right)\cap\M$ is the subset of $\M$ where $X$ can
  re-enter the manifold after having hit the boundary; this is not
  really surprising, since the process behaves like a diffusion
  process\footnote{This is not a diffusion process, strictly
    speaking, since the paths are not continuous.} locally in $\M$.
\item The well-known absorbing boundary condition, $p_t=0$, holds
  on the subset of the boundary where at least one of the vector
  fields $A_r$, $1 \leq r \leq n$, is not tangential.
\item There is a conservation equation for the probability mass
  flowing through the reset map $\Phi$, which is naturally expressed
  using the probability current~$J_t$.
\end{enumerate}
This will be stated more precisely in Theorem~\ref{thm:Main}, which is
our main result. Such a PDE was given by Malham{\'e} et al.
\cite{malhame85electric}, for an example of one-dimensional SHS with
two discrete states. Their proof is partly heuristic, especially
concerning the boundary conditions. Our result extends their work in
several directions, allowing for multi-dimensional domains and
terminal states.

The paper is organized as follows: in section
\ref{section:model_and_propoerties} we give a more precise definition
of SDEs with boundary hitting resets, together with some of their
basic properties.  In the following section, we state and prove our
main result, which gives the FPE equation for solutions of SDEs with
boundary hitting resets. Finally, we illustrate the result with two
applications, the first of which extends Malham{\'e}'s one-dimensional
example.

\section{The model and its basic properties}
\label{section:model_and_propoerties}

\subsection{SDE with boundary hitting resets}

Let $\left( \Omega,\F,(\F_t)_{t\geq0},\Pro{} \right)$ be a filtered space
carrying a $n$-dimensional $(\F_t)$-Brownian motion $B$ and $\left(
  X_t \right)_{t \geq 0}$ a right-continuous $\St$-valued adapted
process.

\pagebreak[5]
\begin{defn}[SDE with boundary hitting resets] \label{defn:sdebhr} We
  say that $X$ solves the SDE with boundary hitting resets
  $\left(\SDE,\Phi\right)$ up to time $\zeta^{-}$ if $\zeta$ is a
  positive stopping time such that, almost surely,
  \begin{enumerate}[(i)\quad]
  \item for all $x\in\T$ and all $s < \zeta$,
    \begin{equation}
      X_s\,=\,x \quad\Rightarrow\quad \forall t \in \left[ s; \zeta \right),\; X_t \,=\, x,
      \jbNoNum
    \end{equation}
  \item $X$ is piecewise continuous on $\left[ 0; \zeta \right) $, i.e. has finitely many
    discontinuities on each compact set $K \subset \left[ 0; \zeta \right) $,
  \item $X$ solves the SDE \eqref{eq:SDE} on each interval of
    continuity $I \subset \left[ 0; \zeta\land\tauEtoile \right)$, where $\tauEtoile$ is the
    first-entrance time of $X$ into $\T$,
  \item and for each discontinuity time $\tau\in [0;\zeta)$, $X$ has a limit on
    the left $X_\tau^-$ such that $X_\tau^- \in\dM$ and $X_\tau = \Phi(X_\tau^-)$.
    \label{property:jumps}
  \end{enumerate}
\end{defn}

\begin{prop}[existence and uniqueness]\label{prop:ExistUnic}
  For each $\F_0$-measurable $\St$-valued random variable $\xi$, there
  is a maximal pathwise unique solution $\left(X,\zeta\right)$ starting
  from $\xi$. More precisely, we mean that:
  \begin{enumerate}[(i)\quad]
  \item $X_0=\xi$ almost surely.
  \item $\zeta$ is a positive stopping time such that $X$ solves
    $\left( \SDE, \Phi \right)$ up to time $\zeta^{-}$.
  \item If $\left( X', \zeta' \right)$ is another solution, then $\zeta' \leq \zeta$
    and $X_t=X_t'$ a.s. on $\left[ 0; \zeta' \right)$.
  \end{enumerate}
\end{prop}

\begin{jbProof}

  \emph{Existence.} We first establish the existence of a solution using
  a recursive construction. More precisely, we construct a sequence of
  right-continuous adapted processes $\left( X^{\indsaut} \right)_{\indsaut \geq
    1}$ and a sequence of stopping times $\left(\tau_{\indsaut}\right)_{\indsaut \geq
    1}$ such that, for all $\indsaut \geq 1$,
  \begin{enumerate}
  \item \label{cond:Existence1} $X_{0}^{\indsaut}=\xi$ almost surely;
  \item $X^{\indsaut}$ has $N_{d}^{\indsaut}$ discontinuities, $0 \leq
    N_{d}^{\indsaut} < \indsaut$, which occur at times $\tau_{1}, \ldots,
    \tau_{N_{d}^{\indsaut}}$;
  \item $X^{\indsaut}$ solves $\left(\SDE,\Phi\right)$ up to time
    $\tau_{\indsaut}^{-}$;
  \item \label{cond:Existence4} if $\tau_{\indsaut} < +\infty$, then
    $X^{\indsaut}$ eventually leaves every compact set of $\M$ before
    time $\tau_{\indsaut}$.
  \end{enumerate}

  Let $Z^{0}$ be the solution\footnote{Thanks to the smoothness of the
    vector fields, equation \eqref{eq:SDE} admits a pathwise unique
    maximal solution on any filtered probability space carrying a
    $n$-dimensional adapted Brownian motion. This solution is defined
    on a random interval $\left[ 0; e \right)$, where $e$ is a
    stopping time called the explosion time.} of \eqref{eq:SDE}
  starting from $\xi$, defined on the event $\left\{ \xi \in \M \right\}$ up
  to the explosion time $e_0$. Exploding solutions are extended to
  $\Rset_{+}$ with the value $\Delta$, where $\Delta \in \T$ is an isolated point
  that acts as a cemetery.  We set
  \begin{equation}
    \tau_{1} \;=\; 
    \left\{
      \begin{array}{lcl}
        e_0 & \quad & \text{if } \xi \in \M, \\
        +\infty & & \text{otherwise},
      \end{array}
    \right.
    \jbNoNum
  \end{equation}
  and
  \begin{equation}
    X_t^{1} \;=\; 
    \left\{
      \begin{array}{lcl}
        Z_t^{0} & \quad & \text{if } \xi \in \M, \\
        \xi & & \text{otherwise}.
      \end{array}
    \right.
    \jbNoNum
  \end{equation}
  It is clear that $X^{1}$ is an adapted process, continuous on
  $\left[ 0; \tau_{1} \right)$, which satisfies the conditions
  \ref{cond:Existence1}--\ref{cond:Existence4}.

  Now we construct $X^{\indsaut+1}$ and $\tau_{\indsaut+1}$ from
  $X^{\indsaut}$ and $\tau_{\indsaut}$, for $\indsaut \geq 1$. We set
  $X^{\indsaut+1}_{t} = X^{\indsaut}_{t}$ on $\left[ 0; \tau_{\indsaut}
  \right)$ and then distinguish between several cases when
  $\tau_{\indsaut} < +\infty$:
  \begin{enumerate}[-]
  \item when $X^{\indsaut}$ has no limit in $\Mb$ at time
    $\tau_{\indsaut}^{-}$, we set $\tau_{\indsaut+1} = \tau_{\indsaut}$ and
    $X^{\indsaut+1}_{t}=\Delta$ for $t \geq \tau_{\indsaut}$;
  \item when $X^{\indsaut}$ has a limit in $\Mb$ at time
    $\tau_{\indsaut}^{-}$ which belongs to $\Phi^{-1} \left( \T \right)$, we
    set $\tau_{\indsaut+1} = +\infty$ and $X^{\indsaut+1}_{t} =
    \Phi(X^{\indsaut}(\tau_{\indsaut}^{-}))$ for $t \geq \tau_{\indsaut}$;
  \item and when $X^{\indsaut}$ has a limit in $\Mb$ at time
    $\tau_{\indsaut}^{-}$ which belongs to $\dM \setminus \Phi^{-1} \left( \M
    \right)$, we set $\tau_{\indsaut+1} = \tau_{\indsaut} + e_{\indsaut}$
    and $X^{\indsaut+1}_{t}=Z^{\indsaut} \left( t - \tau_{\indsaut}
    \right)$, where $Z^{\indsaut}_t$ is the solution of
    $\eqref{eq:SDE}$ with respect to the Brownian motion
    $\tilde{B}^{\indsaut}_{u} = B_{\tau_{\indsaut}+u} -
    B_{\tau_{\indsaut}}$, starting from
    $\Phi(X^{\indsaut}(\tau_{\indsaut}^{-}))$ and defined up to the
    explosion time $e_{\indsaut}$. Note that $\tilde{B}^{\indsaut}$ is
    a Brownian motion with respect to the filtration
    $(\F_{\tau_{\indsaut}+u})_{u\geq 0}$.
  \end{enumerate}
  It is not difficult to check that $X^{\indsaut+1}$ satisfies the
  conditions \ref{cond:Existence1}--\ref{cond:Existence4}. The only
  technical point is to verify that it is a $(\F_t)$-adapted process:
  this comes from the fact that $Z^{\indsaut}$ is progressively
  measurable with respect to the filtration $(\F_{\tau_{\indsaut}+u})_{u\geq
    0}$. Everything else is a direct consequence of our construction.
  Finally, setting $\zeta = \lim_{\indsaut \to +\infty} \tau_{\indsaut}$ and $X_{t}
  = \lim_{\indsaut \to +\infty} X^{\indsaut}(t)$, we obtain a solution of
  $\left(\SDE,\Phi\right)$ up to time $\zeta^{-}$, starting from $\xi$.

  \emph{Pathwise uniqueness.} Let $\left( X',\zeta' \right)$ be another
  solution starting from $\xi$. We use that the solution of
  \eqref{eq:SDE} with a given initial condition is pathwise unique, to
  deduce by recurrence that $X=X'$ almost surely on the interval
  $\left[ 0; \tau_{\indsaut} \land \zeta' \right)$, for all $\indsaut \geq 1$, and
  therefore $X=X'$ almost surely on $\left[ 0; \zeta \land \zeta' \right)$ since
  $\tau_{\indsaut} \to \zeta$. Then we observe that the solution $X$ which has
  been constructed in the first part of the proof cannot be extended
  beyond the lifetime $\zeta$. Indeed, when $\zeta < +\infty$, one of the following
  two situations takes place:
  \begin{enumerate}[-]
  \item if $N_{d} = \lim_{\indsaut \to +\infty} N_{d}^{\indsaut} < +\infty$, the
    solution is exploding in $\M$ at time $\tau_{N_{d}+1}^{-}$ and has no
    limit in $\Mb$;
  \item if $N_{d} = +\infty$, the sequence $\tau_{\indsaut}$ has an
    accumulation point at time $\zeta^{-}$ and therefore cannot be
    extended beyond $\zeta$ without losing the piecewise continuity.
  \end{enumerate}
  This implies that $\zeta' \leq \zeta$, which completes the proof.
\end{jbProof}

\subsection{Basic properties}

Proposition \ref{prop:ExistUnic} shows that, in general, the solution
of a SDE with boundary hitting resets cannot be defined for all $t \geq
0$, either because the solution of the SDE behaves badly or because
the process undergoes an infinite number of resets within a finite
time\footnote{This is what people in the hybrid systems community call
  the Zeno phenomenon \cite{zhang:2001:zeno}, in reference to Zeno of
  Elea and his famous paradoxes.}. Let $\left( X^{x}, \zeta^{x} \right)$
be the maximal solution of $\left( \SDE, \Phi \right)$ starting from $x \in
\St$. We assume from now on that
\begin{assumptions}
\item\label{ass:Conservative} for all $x\in\St$, the solution $X^{x}$ has an infinite lifetime,
  i.e. $\zeta^{x}=+\infty$ almost surely.
\end{assumptions}
Let $( R, \mathcal{R}, \left(\mathcal{R}_t )_{t \geq 0} \right)$ be the
canonical filtered space of all right-continuous paths $\omega:t\in\Rset_{+}
\mapsto \omega(t) \in \St$, with coordinate process $X_{t}:\omega \mapsto \omega(t), t \geq 0$. We set
$\tilde{\Omega}=\Omega\times R$ and $\tilde{\F}=\F\otimes\mathcal{R}$. For all $x \in \St $,
we denote by $\Pro{x}$ the probability law on $(\tilde{\Omega},\tilde{\F})$
defined as the image of $\Pro{}$ by the application
$(\text{Id},X^{x})$. In other words, $\Pro{x}$ is the extension of
$\Pro{}$ such that $X$ is the solution of $( \SDE, \Phi )$ starting from
$x$. The space $( \tilde{\Omega}, \tilde{\F} )$ is endowed with the
filtration $( \tilde{\F}_{t} )_{t \geq 0}$, where $\tilde{\F}_{t} = \F_t
\otimes \mathcal{R}_t$.

Let $\Stb$ be the extended state space obtained as the union of the
state space $\St$ and the boundary $\dM$, i.e. $\Stb \;=\; \St \cup \dM
\;=\; \Mb \cup \T$. We denote by $\Ckc{\dd}{\Mb}$ the set of compactly
supported functions of class $\Cclass^k$ on $\Mb$. For convenience, we
use the notation $\Ck{\dd}{\Stb}$, resp.  $\Ckc{\dd}{\Stb}$, to denote
the set of all functions $\varphi:\Stb\to\Rset$ whose restriction to $\Mb$
belongs to $\Ck{\dd}{\Mb}$, resp.  $\Ckc{\dd}{\Mb}$. Moreover, we
extend the vector fields $\VF_r$, $0\leq r\leq n$, setting
\begin{equation}
  \left(A_r \varphi \right) \left( x \right) = 0,\quad \text{ for all } x\in\T.
\end{equation}

As in the proof of Proposition \ref{prop:ExistUnic}, we denote
$\tau_{\indsaut}$ the stopping time corresponding to the $j^\text{th}$
jump, $\indsaut \geq 1$, with the convention that $\tau_{\indsaut} = +\infty$ if
the process has less than $\indsaut$ jumps.

\pagebreak[5]
\begin{prop}
  The family $\left\{ \Pro{x},\, x\in\St \right\}$ is a Markovian system
  of probability measures on $(\tilde{\Omega}, \tilde{\F},
  (\tilde{\F}_{t})_{t \geq 0})$ with the following properties:
  \begin{enumerate}[(i)]
  \item For all $x\in\St$, $\varphi\in\Ck{2}{\Stb}$ and $t \geq 0$, it holds
    $\Pro{x}$-almost surely that 
    \begin{equation}\label{tatta2}
      \begin{split}
        \varphi(X_t) \;=\; \varphi(x) \,+\, \int_0^t \left( \VF_0 \varphi \right) \left(
          X_s \right)\, ds
        & \,+\, \sum_{r=1}^n\, \int_0^t \left( \VF_r \varphi \right) \left( X_s \right)\, \dBStrato_s^r \\
        & \,+\, \sum_{\tau_\indsaut \leq t}\, \left( \varphi\circ\Phi - \varphi \right) \left(
          X_{\tau_\indsaut}^{-} \right) \,.
      \end{split}
    \end{equation}
  \item Let $\Loi{0}$ be any probability measure on $\St$. Then, for
    all $\varphi\in\Ckc{2}{\Stb}$ and for all $t\geq 0$,
  \begin{equation}
    \label{eq:DynkinFormula}
    \begin{split}
      \Espx{\Loi{0}}{\varphi(X_t)} \;=\; \Loi{0}(\varphi)
      & \,+\, \Espx{\Loi{0}}{\int_0^t \left( \OpKolBack\varphi \right) \left( X_s \right)\, ds} \\
      & \,+\, \Espx{\Loi{0}}{\sum_{\tau_j \leq t}\, \left( \varphi\circ\Phi - \varphi \right) \left(
          X_{\tau_j}^{-} \right)} ,
    \end{split}
  \end{equation}
  where $\OpKolBack$ is given by \eqref{eq:Hormander} and
  $\EspLetter_{\Loi{0}}$ denotes as usual the expectation with respect
  to the probability measure $\Pro{\Loi{0}}=\int_{\St}\, \Loi{0}(\mathrm{d}
  x)\, \Pro{x}$.
  \end{enumerate}
\end{prop}

\begin{jbProof}
  Using the usual chain rule on the intervals of continuity and a
  recurrence on $k$, it is easy to show that, for all $k \geq 1$,
  $\Pro{x}$-almost surely,
  \begin{equation}
    \begin{split}      
      \varphi(X_t) \;=\; \varphi(x) \,+\, \int_0^t \left(\VF_0 \varphi\right)
      \left(X_s\right)\, ds
      & \,+\, \sum_{r=1}^n\, \int_0^t \left( \VF_r \varphi \right) \left( X_s \right)\, \dBStrato_s^r \\
      & \,+\, \sum_{\tau_\indsaut \leq t}\, \left( \varphi \left( X_{\tau_\indsaut}
        \right) - \varphi \left( X_{\tau_\indsaut}^{-} \right) \right)
    \end{split}
    \jbNoNum
  \end{equation}
  on $[0;\tau_{\indsaut})$. This equation holds in fact on $[0;+\infty)$ since
  $\tau_{\indsaut} \to +\infty$. Moreover, we have that $X_{\tau_\indsaut} =
  \Phi(X_{\tau_\indsaut}^-)$ by
  Definition~\ref{defn:sdebhr}(\ref{property:jumps}), which establishes~\eqref{tatta2}. 

  Now we assume that $\varphi$ is compactly supported. The connection
  between Stratonovich integrals and \Ito integrals (see e.g.
  \cite[Theorem V.1.2]{IW89}) yields
  \begin{equation}
    \left( \VF_r \varphi \right) \left( X_s \right)\, \dBStrato_s^r \;=\;
    \left( \VF_r \varphi \right) \left( X_s \right)\, \dDiff B_s^r 
    + \frac{1}{2}\, \left( \VF_r^2 \varphi \right) \left( X_s \right)\, \dDiff s \,,\quad    
    1 \leq r \leq n \,.
    \jbNoNum
  \end{equation}
  Therefore, equation~\eqref{tatta2} becomes
  \begin{equation}
    \begin{split}
      \varphi(X_t) \;=\; \varphi(x) \,+\, \int_0^t \left( \OpKolBack \varphi \right) \left( X_s
      \right)\, ds
      & \,+\, \sum_{r=1}^n\, \int_0^t \left( \VF_r \varphi \right) \left( X_s \right)\, \dDiff B_s^r \\
      & \,+\, \sum_{\tau_\indsaut \leq t}\, \left( \varphi\circ\Phi - \varphi \right) \left(
        X_{\tau_\indsaut}^{-} \right) \,,
    \end{split}
    \jbNoNum
  \end{equation}
  and equation~\eqref{eq:DynkinFormula} is obtained by taking
  expectations and using that the process $M_t^{\varphi} = \sum_{1\leq r \leq n}
  \int_{0}^{t} (\VF_r \varphi) (X_s) \dDiff B_s^r$ is a martingale (since
  $\VF_r \varphi$ is bounded).

  \emph{Proof of the Markov property.} We define the semigroup
  $(P_t)_{t\geq 0}$, as usual, by $(P_t \varphi)(x)=\Espx{x}{\varphi(X_t)}$. Then,
  for all $s<t$, using equation~\eqref{tatta2} and the martingale
  property of $M_t^{\varphi}$,
  \begin{align}
    &\; \Espx{x}{\varphi(X_t) \,|\, \tilde{\F}_s} \\
    =\;&\; \EspLetter_{x} \left\{ \varphi(X_s) \,+\, \int_s^t \left( \OpKolBack\varphi \right) \left( X_s \right)\, ds \,+\,
      \sum_{s<\tau_j \leq t}\, \left( \varphi\circ\Phi - \varphi \right) \left( X_{\tau_j}^{-} \right) \,|\, \tilde{\F}_s \right\} \jbNoNum\\
    =\; & (P_{t-s}\varphi)(X_s) \,.\jbNoNum
  \end{align}
  The last line follows from \eqref{eq:DynkinFormula} and the fact
  that, under the probability $\Pro{x} \{\, \cdot \,|\, \tilde{\F}_s \}$,
  the process $(X_{s+u})_{u\geq0}$ is the maximal solution of $(\SDE,\Phi)$
  starting from $X_s$ and driven by $(B_{s+u})_{u\geq0}$.  Therefore, it
  has the same distribution as $(X,B)$ under $\Pro{X_s}$ (this is a
  consequence of the uniqueness in law for \eqref{eq:SDE}, using a
  proof similar to that of Proposition~\ref{prop:ExistUnic}).
\end{jbProof}

\begin{jbRemark} The same Markov process can also be obtained using
  the ``revival theorem'' of Ikeda-Nagasawa-Watanabe
  \cite{INW66,Mey75}. This can be found with full details in
  \cite{Bujorianu:2004:TheoreticalFoundations}, in a framework which
  is much more general than ours. The idea is to restart the process
  $X^{1}$ in the proof of Proposition~\ref{prop:ExistUnic}, according
  to the deterministic kernel $N(\omega,.) \;=\; \delta_{\Phi(X^{1}(\tau_{1}^{-}(\omega),
    \omega))}$. We have preferred a more explicit pathwise construction,
  which, we hope, provides a better insight into the behaviour of the
  process. 
\end{jbRemark}

\section{Fokker-Planck-Kolmogorov's forward equation}
\label{section:FPKfe}

\subsection{Main result of the paper}

Let $\Loi{0}$ be a probability measure on $\St$, and denote
by $\mu_t$ the probability law of $X_t$ at time $t\geq 0$ under $\Pro{\Loi{0}}$. We
assume that
\begin{assumptions}
\item\label{Ass:pdf_exists} for all $t\geq 0$, $\Loi{t}$ has a density
  $p_t=p(\cdot,t)$ on $\M$, with respect to the Riemannian volume measure $\rvm$. 
\end{assumptions}
Consequently, we can decompose $\mu_t$ as
\begin{equation}
  \label{eq:Law}
  \mu_t \;=\; p_{t}\, \rvm \;+\; \sum_{x\in\T} q(x,t)\, \delta_x , 
\end{equation}
where $q(x,t) \;=\; \mu_t\left( \{x\} \right) \;=\; \Pro{\Loi{0}} \left\{
  X_t=x \right\}$.

\medskip\emph{Assumptions about the boundary and the reset function.}
Set $\Gua = \Phi^{-1} \left( \Hyp\right)$. Since the boundary is
zero-dimensional when $n=1$, our assumptions will be slightly
different in the one-dimensional case and in the multi-dimensional
case. First, regardless of the dimension, we assume that
\begin{assumptions}
\item\label{ass:Bijection} $\Phi_{|\Gua}$ is a bijection from $G$ to
  $H$\footnote{This assumption could easily be relaxed in our proof,
    cf. Remark \ref{rem:PostDemo}.}.
\end{assumptions}
In the multi-dimensional case, we further assume that
\begin{assumptions}
\item\label{ass:Reset} $\Hyp$ is a closed and orientable
  hypersurface\footnote{embedded smooth submanifold, without boundary
    and of codimension 1}, and $\Phi_{|\Gua}$ is a
  \mbox{$\Cclass^2$--diffeomorphism} from $\Gua$ to $\Hyp$.
\end{assumptions}
Let $\sigma$ be the surface measure induced by $\scalfct$ on $\dM$ and
$\Hyp$, when $n \geq 2$. We define $h=\left|\jacob \Phi\right|$ on $\Gua$,
where $\jacob \Phi$ is the Jacobian of $\Phi$ with respect to the Riemannian
volume forms on $\Gua$ and $\Hyp$ (for any choice of orientation),
such that, for all $f \in \Ckc{0}{\Hyp}$,
\begin{equation}
  \label{eq:ChangeOfVariable}
  \int_{\Hyp}\, f \dsigma \;=\; \int_{\Gua}\, f\circ\Phi\; h \dsigma
\end{equation}
The same formula holds in the one-dimensional case if we set $h=1$ and
interpret $\sigma$ as the counting measure on $\dM \cup \Hyp$.

\emph{Smoothness assumptions for the probability law.} Finally, we
assume that $p$ and $q$ are smooth enough, that is:
\begin{assumptions}
\item\label{ass:RegDensity} $p$ is of class $\Cclass^{2,1}$ on $\left(
    \MHyp \right) \times \Rset_+$. Morevover, for all $t \geq 0$, $p_t$ is
  continuous on $\Mb$, and the differential $\dpdf_t$ is continuous on
  $\Mb \setminus \Hyp$ with a discontinuity of the first kind on $\Hyp$.
\item\label{ass:Reg_q} For each $x\in\T$, $t \mapsto q\left(x,t\right)$ is
  continuously differentiable on $\Rset_+$.
\end{assumptions}
Note that we have not assumed $\dpdf_t$ continuous on $\Mb$, since it
is precisely the discontinuities of $\dpdf_t$ through $H$ that account
for the jumping behaviour of $X$. Assumption \ref{ass:RegDensity}
implies that, for all $t \geq 0$, the probability current $J_t$ defined
by \eqref{eq:Jdef} is a $\Cclass^{1}$ vector field on $\MHyp$, that has
a limit on $\dM$ and a discontinuity of the first kind on $\Hyp$. We
set $\Jout_{t} \; = \; \scal{J_t}{\nu}$ on $\dM$, where $\nu$ is the
outward-pointing unit normal, and
\begin{equation}
  \begin{split}
    \Jin_{t} & = \langle\, J_{t}^{(1)} ,\, \nu_{21} \,\rangle \,+\,  \langle\, J_{t}^{(2)} ,\, \nu_{12} \,\rangle  \\
    & = \langle\, J_{t}^{(2)} - J_{t}^{(1)} ,\, \nu_{12} \,\rangle
  \end{split}
\end{equation}
on $\Hyp$, where $\nu_{12}=-\nu_{21}$ is the unit normal to $\Hyp$
directed from side $1$ to side $2$, and $J_{t}^{(i)}$ is the value of the
discontinuous vector field $J_{t}$ on the side $i$ of $H$, $i\in\{1,2\}$. The
first expression makes it clear that the definition of $\Jin_{t}$ does not
depend on the choice of an orientation\footnote{In fact, $\Hyp$ does
  not even need to be orientable for this to be defined
  \cite{Schwartz:1966:Distributions}.} for $\Hyp$.

\pagebreak[5]
\begin{thm}
\label{thm:Main}
Under assumptions \ref{ass:Conservative}--\ref{ass:Reg_q}, the law
$\mu_t$ of $X$ evolves according to the following equations:
\begin{align}
  \frac{\partial p}{\partial t}(x,t) \;&=\; (\OpKolForw p_t)(x)  & & \text{on } \MHyp \,,
  \label{eq:FPKE} \\
  \frac{\partial q}{\partial t}(x,t) \;&=\; \int_{\Gua_x} \Jout_{t}\, \dDiff \sigma & &
  \text{on } \T \,,
  \label{eq:Conserv1}
\end{align}
where $\Gua_x=\Phi^{-1} \left( \left\{ x \right\} \right)$, for all $x \in \T$.  Moreover,
the following boundary conditions hold:
\begin{align}
  \Jout_{t} \;&=\; h\; \Jin_{t} \circ \Phi\, & & \text{on } G \,,
  \label{eq:Conserv2} \\
  p_{t} \;&=\; 0 & & \text{on } \dM_0 \,,
  \label{eq:BoundCond}
\end{align}
where $\dM_{0}$ is the subset of the boundary defined by
\begin{equation}
  \dM_{0} \;=\; \left\{\, x\in\dM \,|\, \exists r \in \left\{ 1, \ldots, n \right\},\;
    \scal{\VF_r}{\nu} \neq 0 \,\right\} \,.
  \label{eq:NonCharBound}
\end{equation}
\end{thm}

\begin{jbRemark}
  $\dM_0$ can be seen as the non-characteristic part of the boundary
  $\dM$ with respect to the operator $\OpKolForw$.  Indeed,
  $\OpKolForw$ is a second order operator with principal symbol
  \begin{align}
    \sigma_2\left(\OpKolForw\right)\left(x,\xi\right) 
    & \;=\; -\frac{1}{2}\, \sum_{r=1}^{n}\, \sum_{i,j=1}^{n}\,
    \VF_{r}^{i}\left(x\right) \VF_{r}^{j}\left(x\right) \xi_i \xi_j 
    \jbNoNum\\
    & \;=\; -\frac{1}{2}\, \sum_{r=1}^{n}\, \left( 
      \sum_{i=1}^{n}\,\VF_{r}^{i}\left(x\right) \xi_i \right)^{2},
    \jbNoNum
  \end{align}
  where the $\xi_i$'s are the coordinates of the covector $\xi$ in any
  local orthonormal coframe. Therefore, $\dM$ is characteristic at $x$
  if and only if $A_r$ is tangent to $\dM$ at $x$, for each $r \in
  \left\{ 1, \ldots, n \right\}$.
\end{jbRemark}

\subsection{Proof of Theorem \ref{thm:Main}}

The result is obvious when $\Loi{0}(\M)=0$, since $X$ is then a
constant process. When $\Loi{0}(\M)>0$, we observe that $\Pro{\Loi{0}}
\{\, \cdot \,|\, X_0\in\M \,\} = \Pro{\Loi{0}'}$, where $\Loi{0}'=\Loi{0} /
\Loi{0}(\M)$. Thus, it can be assumed without loss of generality that
$\Loi{0}(\M)=1$, i.e. that $q(x,0)=0$ for all $x\in\T$.

$\triangleright\;$ As a first step, we will prove that, for each test function
$\varphi\in\Ckc{2}{\Stb}$ such that
\begin{equation}
  \label{eq:CondComp}
  \varphi\circ\Phi=\varphi \quad \text{ on } \Gua , 
\end{equation}
the following equation holds:
\begin{equation}
  \label{eq:FirstPartOfTheProof}
  \int_{\M} \varphi\, \left( p_t - p_0 \right) \dm 
  \; = \; \int_{0}^{t} \int_{\M} L\varphi\, p_s \dm\; \ds \; - \; 
  \Espx{\Loi{0}}{ 1_{\left(\tauEtoile\leq t\right)}\,
    \varphi\left(X_{\tauEtoile}^{-}\right) } \,.
\end{equation}
Indeed, \eqref{eq:CondComp} implies that the jump term $\left( \varphi\circ\Phi - \varphi
\right) (X_{\tau_k}^-)$ in equation \eqref{eq:DynkinFormula} vanishes
when the process undergoes a reset into $\M$ at time $\tau_k$.  Hence,
there is at most one term left in the summation, corresponding to a
possible jump from $\dM$ to a terminal state in $\T$ at time
$\tauEtoile$:
\begin{align}
  \Espx{\Loi{0}}{\sum_{\tau_\indsaut \leq t}\, \left( \varphi\circ\Phi - \varphi \right) \left(
      X_{\tau_\indsaut}^- \right)} &\;=\; \Espx{\Loi{0}}{
    1_{\left(\tauEtoile\leq t\right)} \left( \varphi\circ\Phi - \varphi \right)
    \left(X_{\tauEtoile}^{-}\right) }
  \jbNoNum\\
  &\;=\; \sum_{x\in\T}\, \varphi(x)\, q(x,t) \,-\, \Espx{\Loi{0}}{
    1_{\left(\tauEtoile\leq t\right)} \varphi \left(X_{\tauEtoile}^{-}\right) }
  \,.
  \label{eq:Term1}
\end{align}
Using Fubini's theorem together with the decomposition \eqref{eq:Law}
of $\mu_t$ and the fact that $\OpKolBack\varphi$ vanishes on $\T$, we can
rewrite the second expectation in \eqref{eq:DynkinFormula} as
\begin{align}
  \Espx{\Loi{0}}{\int_0^t \left( \OpKolBack\varphi \right) \left( X_s \right)\,
    \ds} & \;=\; \iint_{\Omega\times \left[ 0;t \right]}\, \left(L\varphi\right)
  \left(X(s,\omega)\right)\; \Pro{\Loi{0}}(\dDiff \omega)\, \ds
  \jbNoNum\\
  & \;=\; \int_{0}^{t} \int_{\St} L\varphi\; \dDiff \mu_s\, \ds
  \jbNoNum\\
  & \;=\; \int_{0}^{t} \int_{\M} L\varphi\; p_s\, \dm\, \ds \,.
  \label{eq:Term2}
\end{align}
Using equation \eqref{eq:Law} once more, we can also expand the
left-hand side of \eqref{eq:DynkinFormula}:
\begin{equation}
  \Espx{\Loi{0}}{\varphi(X_t)} \;=\; \sum_{x\in\T}\, \varphi(x)\, q(x,t)
  \;+\; \int_{\M}\, \varphi\, p_t\, \dm \,.
  \label{eq:Term3}
\end{equation}
Finally, replacing expressions \eqref{eq:Term1}, \eqref{eq:Term2},
and \eqref{eq:Term3} back into \eqref{eq:DynkinFormula} completes the
proof of equation \eqref{eq:FirstPartOfTheProof}. 

$\triangleright\;$ The next step of the proof relies on the following version of
Stokes's formula:
\begin{equation}
  \label{eq:Stokes}
  \int_{\M\setminus\Hyp} \diverg A\, \dm \;=\; \int_{\dM} \scal{A}{\nu}\, \dsigma
  \;-\; \int_{\Hyp} \scal{A^{(2)} - A^{(1)}}{\nu_{12}}\, \dsigma \,,
\end{equation}
where $A$ is a compactly supported vector field on $\Mb$, of class
$\Cclass^1$ in $\MHyp$, extending continuously to $\dM$ and having
a discontinuity of the first kind on $\Hyp$. This formula is easily
proved using the usual divergence theorem \cite[Theorem
14.34]{Lee:2003:SmoothManifolds} or the general results of
\cite[Chapter IX]{Schwartz:1966:Distributions}. Using this formula
together with the assumption \ref{ass:RegDensity} that $p_t$ is
continuous on $\Mb$, for all $t \geq 0$, yields:
\begin{equation}
  \label{eq:ResultatIPP}
  \begin{split}
    \int_{\M}\, \OpKolBack \varphi\, p_t\, \dm \; = \; & \int_{\M\setminus\Hyp}\, \varphi\,
    \OpKolForw p_t\, \dm \; + \; \beta_{t}(\varphi) \,,
  \end{split}
\end{equation}
where $\OpKolForw$ is the formal adjoint of $\OpKolBack$, i.e.
\begin{equation}
  \OpKolForw p_t \;=\; -\, \diverg\left(p_t\VF_0\right)
  \;+\; \frac{1}{2}\, \sum_{r=1}^{n}\, \diverg\left(\diverg\left(p_t\VF_r\right)A_r\right) \,,
  \jbNoNum
\end{equation}
and $\beta_{t}$ is a distribution supported by $\dM\cup\Hyp$:
\begin{equation}
  \label{eq:BoundaryDistribution}
  \beta_{t}(\varphi) \;=\; \frac{1}{2}\, \sum_{r=1}^{n}\, \int_{\dM}\,
  \VF_r \varphi\; p_t\, \scal{\VF_{r}}{\nu}\, \dsigma \;+\; \int_{\dM}\, \varphi\, \Jout_t\, \dsigma - \int_{H}\, \varphi\, \Jin_t\, \dsigma \,.
\end{equation}
A detailed derivation of these equations can be found in
appendix~\ref{appendix:IPP}. Using \eqref{eq:ResultatIPP}, we rewrite
equation~\eqref{eq:FirstPartOfTheProof} as
\begin{equation}
  \label{eq:SecondPartOfTheProof}
  \int_{\M\setminus\Hyp}\, \varphi\, \tilde{p}_{t}\; \dm  \;=\; \int_{0}^{t}\, \beta_{s}(\varphi)\, \ds \;-\; \Espx{\Loi{0}}{ 1_{\left(\tauEtoile\leq t\right)}\,
    \varphi\left(X_{\tauEtoile}^{-}\right) } \,,
\end{equation}
where
\begin{equation}
  \tilde{p}_{t} \;=\; p_t - p_0 - \int_{0}^{t}\, \OpKolForw p_s\, \ds
  \jbNoNum
\end{equation}
is well-defined and continuous on $\MHyp$.

$\triangleright\;$ The last step of the proof consists in choosing specific test
functions $\varphi$ in equation \eqref{eq:SecondPartOfTheProof}. First, it
is easy to check that, if $\varphi$ has its support in $\MHyp$, then
\eqref{eq:CondComp}~is automatically satisfied since $\varphi$
vanishes on both $\Gua$ and $\Hyp$, and equation
\eqref{eq:SecondPartOfTheProof} becomes
\begin{equation}
  \int_{\M\setminus\Hyp}\, \varphi\,\tilde{p}_t\; \dm  \;=\; 0,
  \jbNoNum
\end{equation}
which proves that $\tilde{p}_t=0$, for all $t\geq0$. Equation
\eqref{eq:FPKE} then follows by differentiation.

Before proceeding to the derivation of the three other equations of
the theorem, we shall state a useful technical lemma, whose proof can
be found in appendix~\ref{appendix:ExtLemma}.

\begin{lem}\label{lem:ExtensionLemma}
  For any $\eta_i \in \Ckc{2}{\dM}$, $i \in \{ 1,2 \}$, such that
  ${\eta_1}_{|\Gua}$ is compactly supported in $\Gua$, we can find a
  function $\varphi \in \Ckc{2}{\Stb}$ such that
  \begin{enumerate}[(i)]
  \item\label{it:L1} $\varphi$ satisfies \eqref{eq:CondComp},
  \item\label{it:L2} $\varphi=\eta_1$ on $\dM$,
  \item\label{it:L3} and $\frac{\partial \varphi}{\partial \nu}=\eta_2$ on $\dM$.
  \end{enumerate}
\end{lem}

\emph{Proof of equation~\eqref{eq:BoundCond}.} For any $\eta \in
\Ckc{2}{\dM}$, Lemma~\ref{lem:ExtensionLemma} provides us with a
function $\varphi \in \Ckc{2}{\Stb}$ such that \eqref{eq:CondComp} is
satisfied, $\varphi=0$ on $\dM \cup \Hyp$ and $\frac{\partial\varphi}{\partial\nu}=\eta$ on $\dM$.  Then
equation~\eqref{eq:SecondPartOfTheProof} implies that $\beta_{t}(\varphi)=0$,
for all $t \geq 0$. Using equation~\eqref{eq:BoundaryDistribution} and
the fact that $\varphi$ vanishes on $\dM \cup \Hyp$, this can be written as
\begin{equation}
  \label{eq:DudulePower}
  \sum_{r=1}^{n}\, \int_{\dM}\, \VF_r \varphi\; p_t\, \scal{\VF_{r}}{\nu}\, \dsigma \;=\; 0.
\end{equation}
Moreover, by construction, $\varphi$ has the property that $\VF_r \varphi \;=\;
\scal{\VF_{r}}{\nu}\; \eta$, for all $r \in \{1, \ldots, n\}$, which allows to rewrite equation~\eqref{eq:DudulePower} as
\begin{equation}
  \int_{\dM}\, \left( \sum_{r=1}^{n}\, \scal{\VF_{r}}{\nu}^{2}\, \right)\, p_t\; \eta\; \dsigma \;=\; 0 \,.
  \jbNoNum
\end{equation}
This holds for all $\eta \in \Ckc{2}{\dM}$, which proves that
$p_t\, \left( \sum_{r=1}^{n}\, \scal{\VF_{r}}{\nu}^{2}\, \right)=0$ on
$\dM$, for all $t\geq 0$. Observing that $\sum_{r=1}^{n}\,
\scal{\VF_{r}}{\nu}^{2} > 0$ on $\dM_{0}$ establishes
equation~\eqref{eq:BoundCond}.

\emph{Proof of equation~\eqref{eq:Conserv2}.} For any $\eta \in
\Ckc{2}{G}$, we can find by Lemma~\ref{lem:ExtensionLemma} a function
$\varphi \in \Ckc{2}{\Stb}$ satisfying~\eqref{eq:CondComp} and
\begin{equation}
  \varphi_{|\dM} \;=\;
  \left\{ 
    \begin{array}{lcl}
      \eta & \quad & \text{on } \Gua, \\
      0 & \quad & \text{on } \dM \setminus \Gua.
    \end{array}
  \right.
  \jbNoNum
\end{equation}
For such a function, equation~\eqref{eq:SecondPartOfTheProof} reduces to
\begin{equation}
  \int_{\Gua}\, \varphi \Jout_t\, \dsigma - \int_{H}\, \varphi\, \Jin_t\, \dsigma = 0 \,.
  \jbNoNum
\end{equation}
Then, since $\varphi$ satisfies \eqref{eq:CondComp} and $\varphi=\eta$ on $\Gua$, we have
\begin{equation}
  \int_{\Gua}\, \eta\, \left( \Jout_t - h\, \Jin_t \circ \Phi \right)\, \dsigma = 0 \,.
  \jbNoNum
\end{equation}
This holds for all $\eta \in \Ckc{2}{G}$, which proves equation~\eqref{eq:Conserv2}.

\emph{Proof of equation~\eqref{eq:Conserv1}.} This time we choose $\eta \in
\Ckc{2}{\Gua_{x}}$ for some $x\in\T$, $\varphi=\eta$ on $\Gua_{x}$ and $\varphi=0$ on
$\dM \setminus \Gua_{x}$. Equation~\eqref{eq:SecondPartOfTheProof} then
becomes
\begin{equation}
  \Espx{\Loi{0}}{
    1_{\left\{ \tauEtoile\leq t,\, X_{\tauEtoile}^{-} \in \Gua_{x} \right\} }\;
    \eta \left(X_{\tauEtoile}^{-}\right)
  } \;=\; \int_{0}^{t} \int_{G_x}\, \Jout_s\, \eta\; \dsigma \ds \,.
  \jbNoNum
\end{equation}
This implies by dominated convergence that
\begin{equation}
  \label{eq:DomConvResult}
  \Pro{\Loi{0}} \left\{ \tauEtoile\leq t,\; X_{\tauEtoile}^{-} \in K \right\} \;=\; \int_{0}^{t} \int_{K}\, \Jout_s\, \dsigma \ds \,,
\end{equation}
for any compact subset $K$ of $\Gua_x$. Since
the left-hand side is increasing with $t$, this shows that $\Jout_s \geq
0$ on $\Gua_x$, for all $s\geq 0$. Consequently, letting $K \uparrow \Gua_{x}$,
equation~\eqref{eq:DomConvResult} yields
\begin{equation}
  \Pro{\Loi{0}} \left\{ \tauEtoile\leq t,\; X_{\tauEtoile}^{-} \in \Gua_x \right\}
  \;=\; \int_{0}^{t} \int_{\Gua_x}\, \Jout_s\, \dsigma \ds
  \jbNoNum
\end{equation}
by monotone convergence. Finally, observing that
\begin{equation}
  \Pro{\Loi{0}} \left\{ \tauEtoile\leq t,\; X_{\tauEtoile}^{-} \in \Gua_x \right\} \;=\; q(x,t)
  \jbNoNum
\end{equation}
yields equation~\eqref{eq:Conserv1} and thus completes the proof of
~Theorem~\ref{thm:Main}.

\begin{jbRemark}
  \label{rem:PostDemo}
  This proof can easily be generalized to the case where $\Phi$ is no
  longer a diffeomorphism, but still a local diffeomorphism such that
  $\Phi^{-1}\left(\{x\}\right)$ is finite for all $x\in\Hyp$. In this case,
  we have
  \begin{equation}
    \Jin_t (x) \;=\; \sum_{y\in\Phi^{-1}\left(\{x\}\right)}\, h^{-1}(y)\, \Jout_t(y)
    \,,\qquad \forall x\in\Hyp,\; \forall t\geq 0,
    \jbNoNum
  \end{equation}
  instead of equation~\eqref{eq:Conserv2}.
\end{jbRemark}

\section{Examples}

\subsection{Stochastic hybrid systems}

As a first example, we will give the FPE for a multi-dimensional SHS
with two discrete states, which models the temperature in a house with
$n$ rooms, $n \geq 1$, regulated by a single thermostat. This is a
generalisation of the one-dimensional process that was studied in
\cite{malhame85electric}.

Let the $n$-tuple $\theta = ( \theta_{1}, \ldots, \theta_{n} ) \in \Rset^{n}$ describe the
temperature in the $n$ rooms of the house, and $q \in \mathcal{Q}= \{ 0,
1 \}$ the binary state of the thermostat. The global state of the
system is then described by the variable $x = (q,\theta) \in
\mathcal{Q}\times\Rset^{n}$, which has both a discrete and a continuous
component, hence the name ``hybrid''. For a given state $q \in
\mathcal{Q}$ of the thermostat, the temperature $\Theta_t$ evolves in
$\Rset^{n}$ according to a SDE of the form 
\begin{equation}
  \tag{$\SDE_{q}$}\label{eq:SDEq}
  \mathrm{d}\Theta_{t} \;=\; f(q,\Theta_t) \dt \;+\; \gamma \dBIto_{t} \,,
\end{equation}
where $\gamma \in \Rset^{n\times n}$ and $f_{q} = f(q,\cdot)$ describes the action of
the thermostat, the effect of the exterior environment, and the
coupling between the temperatures of adjacent rooms. The switching of
the thermostat is controlled by a linear criterion $\Psi(\theta)=\sum_{i=1}^{n}
\alpha_{i} \theta_{i}$. The thermostat switches on when $\Psi(\Theta_t)$ crosses some
threshold $\Psi_{\mathrm{min}}$ downwards, and switches off when it
crosses another threshold $\Psi_{\mathrm{max}} > \Psi_{\mathrm{min}}$
upwards. This can be described in the SHS framework as follows: we define
\begin{align}
  M_0 & \;=\; \left\{ \, \theta\in\Rset^{n} \,|\, \Psi(\theta) > \Psi_{\mathrm{min}}
    \,\right\} \,, \jbNoNum\\
  M_1 & \;=\; \left\{\, \theta\in\Rset^{n} \,|\, \Psi(\theta) < \Psi_{\mathrm{max}}
    \,\right\} \,,\jbNoNum\\
  \intertext{and then the so-called ``hybrid state space''} \Mb &
  \;=\; \{0\}\times\overline{M}_0 \;\cup\; \{1\}\times\overline{M}_1 \;\subset\; \mathcal{Q}
  \times \Rset^{n} \,. \jbNoNum
\end{align}
The process $\Theta_t$ evolves continuously in $M_q$ according to the SDE
$\SDE_q$ as long as the thermostat is in state $Q_t=q$, and $Q_t$
switches when $\Theta_t$ reaches $\partial M_q$. Therefore, the hybrid process
$X_t = (Q_t,\Theta_t)$ takes its values in $\M$. Considering the atlas
\begin{equation}
  \mathcal{A} \;=\; \left\{ (\{0\}\times\overline{M}_0, \pi), (\{1\}\times\overline{M}_1, \pi) \right\} \,,
  \jbNoNum
\end{equation}
where we denote by $\pi$ the projection $\pi:(q,\theta)\mapsto \theta$, $\Mb$ can be seen
as a smooth manifold with boundary, which has two components. The
process $X_t$ is then the solution of a SDE with boundary hitting
resets $(\SDE,\Phi)$, with the reset map
\begin{align}
  \Phi:\quad \dM &\;\longrightarrow \; \M \jbNoNum\\
  (q,\theta) &\; \longmapsto \; (1-q,\theta) \jbNoNum
\end{align}
and the vector fields
\begin{align}
  \VF_0(x) &\;=\; \sum_{i=1}^{n}\, f^{i}(x)\, \frac{\partial}{\partial\theta_i} \,, \jbNoNum\\
  \VF_r(x) &\;=\; \sum_{i=1}^{n}\, \gamma^{i}_{r}\, \frac{\partial}{\partial\theta_i} \,,\quad 1
  \leq r \leq n \,. \jbNoNum
\end{align}
We note that it is not always possible to view a hybrid state space
$\M = \cup_{q\in\mathcal{Q}} \{q\}\times M_q$ as a manifold, since the domains
$M_q$ have possibly different dimensions in the general SHS framework
\cite{Hu:2000:TTSHS}.

We will now give the FPE equation associated to $(\SDE,\Phi)$, assuming a
priori that the assumptions \ref{ass:Conservative},
\ref{Ass:pdf_exists} and \ref{ass:RegDensity} are satisfied
(\ref{ass:Bijection} and \ref{ass:Reset} are easily checked, and
\ref{ass:Reg_q} is meaningless here since there is no terminal state,
i.e. $\T=\varnothing$). The metric on $\Mb$ is the one induced by the
Euclidian metric on each $M_q$, $q\in\mathcal{Q}$. The density $p_t$ can
be seen as a pair $(p_t^0, p_t^1)$ of modal densities, where $p_t^q$
is defined on $\overline{M}_q$. These modal densities solve a system
of PDE, 
\begin{equation}
  \frac{\partial p_t^q}{\partial t} \;=\; -\sum_{i=1}^{n}\, \frac{\partial f^i_q p_t^q}{\partial\theta_i}
  \;+\; \frac{1}{2}\, \sum_{i,j=1}^{n}\, a^{ij}\, \frac{\partial^2 p_t^q}{\partial\theta_i \partial\theta_j} 
  \,,\quad q\in\mathcal{Q} \,,
  \jbNoNum
\end{equation}
where $a^{ij}=\sum_{r=1}^{n} \gamma_r^i \gamma_r^j$. Assuming that the matrix $a$
is positive, each modal density satisfies the absorbing boundary
condition $p_t^q=0$ on $\partial M_q$. Moreover, the PDEs are coupled by the
conservation equations
\begin{equation}
  \label{eq:ConservThermo}
  \Jout_t (q,\theta) \;=\; \Jin_t (1-q,\theta) \,, \qquad q\in\mathcal{Q} ,\, \theta \in \partial M_q \,,
\end{equation}
where the probability current $J_t=\sum_{i=1}^{n} J_t^i \frac{\partial}{\partial\theta_i}$ is given by 
\begin{equation}
  J_t^i (q,\theta) \;=\; f^i(q,\theta)\, p_t(q,\theta) \;-\; \frac{1}{2}\, \sum_{j=1}^{n}\, a^{ij}\, \frac{\partial p_t}{\partial \theta_j}(q,\theta) \,.
  \jbNoNum
\end{equation}
As in \cite{malhame85electric}, we observe that the drift $f$ does not
appear in the conservation equations \eqref{eq:ConservThermo}, since
the absorbing boundary condition holds.

\subsection{A first exit problem}

Our second application deals with the following problem: let
$X=(X_t)_{t\geq 0}$ be the solution of a SDE of the form \eqref{eq:SDE}
on $\Rset^n$; given an open subset $U\subset\Rset^n$ and an initial
probability law $\Loi{0}$ such that $\supp \Loi{0} \subset U$, we want to
compute
\begin{equation}
  r_i(t) \;=\; \Pro{\Loi{0}} \left\{ \tau \leq t, X_{\tau} \in \partial U_i \right\},
  \qquad 1\leq i\leq D,\quad 0\leq t\leq T
  \jbNoNum
\end{equation}
where $\tau$ is the first exit time of $X$ from $U$, $D\in\Nset^*$, $T>0$,
and $\{ \partial U_i,\, 1\leq i\leq D\}$ is a partition of the boundary $\partial U$. We
assume that the closure $\overline{U}$ of $U$ in $\Rset^n$ is a smooth
manifold with boundary, whose interior coincides with $U$.

It is well-known \cite[Section 5.4]{Schuss:1980:TASDE} that, if the PDE
\begin{align}
  \frac{\partial u_i}{\partial t} & \;=\; \OpKolBack u_i && \text{on } U\times[0;T], \jbNoNum\\
  u_i & \;=\; 0 && \text{on } U\times\{0\}, \jbNoNum\\
  u_i & \;=\; 1_{\partial U_i} && \text{on } \partial U\times(0;T], \jbNoNum
\end{align}
has a bounded solution, then
\begin{equation}
  u_i (x,t) \;=\; \Pro{x} \left\{  \tau \leq t, X_{\tau} \in \partial U_i \right\} \,.
  \jbNoNum  
\end{equation}
This provides a first approach to our problem, since $r_i$ can be
recovered from $u_i$ using an integration with respect to $\Loi{0}$.
Another possible approach is provided by Theorem~\ref{thm:Main}: we
introduce a set of isolated terminal states $\T=\{1, \ldots, D\}$ and consider
the process $\widetilde{X}$ which coincides with $X$ up to time $\tau^-$
and then goes to the state $i\in\T$ such that $X_{\tau}^- \in \partial U_i$. This new
process $\widetilde{X}$ is the solution of a SDE with boundary hitting
resets $(\SDE,\Phi)$ on $U\cup\T$, where the reset map $\Phi$ is defined by
$\Phi(x)=i$ for all $x\in\partial U_i$, $1\leq i\leq D$. The functions $r_i$ can be interpreted in
this framework as
\begin{equation}
  r_i(t) \;=\; \Pro{\Loi{0}} \left\{ \widetilde{X}_{t} = i \right\} \;=\; q(i,t)\,.
  \jbNoNum
\end{equation}
Assuming that the density exists and is smooth enough, the $r_i$'s can be
obtained according to Theorem~\ref{thm:Main} by solving
\begin{align}
  \frac{\partial p}{\partial t} &\;=\; \OpKolForw p && \text{on } U\times[0;T], \jbNoNum\\
  p &\;=\;0 && \text{on } \partial U^{*}\times[0;T], \jbNoNum\\
  \frac{\mathrm{d} r_i}{\mathrm{d}t} &\;=\; \int_{\partial U_i}\, \Jout_t\,
  \dsigma && \text{on } [0;T]. \jbNoNum
\end{align}
where $\partial U^{*}$ is defined as $\dM_{0}$ in Theorem~\ref{thm:Main}. We
point out that our approach based on the forward equation requires the
resolution of a single PDE, in contrast with the first approach, based
on the backward equation, which involves the resolution of $D$ PDEs.
The drawback is that, using the forward approach, the solution is only
obtained for the given initial distribution $\Loi{0}$.

\begin{jbRemark}
  \begin{enumerate}[a)\ ]
  \item From a practical point of view, such a method is of course
    limited to processes with a state space of low dimension, where
    the numerical resolution of the PDEs is feasible. 
  \item A similar methodology can be used to tackle the problem of
    reachability analysis for stochastic hybrid systems (see
    \cite{Bujorianu:2003:DSHS1} for further details), at least when
    the target set is a closed subset with smooth boundary.
  \end{enumerate}
\end{jbRemark}

\section{Appendices}
\subsection{Detailed derivation of equation~\eqref{eq:ResultatIPP}}
\label{appendix:IPP}

In this appendix we give a more detailed derivation of equation
\eqref{eq:ResultatIPP}. First, we choose $\VF = \VF_r \varphi\, p_{t} \VF_r$
in \eqref{eq:Stokes}. Using that $\VF$ is continuous on $\Mb$ and that
$\diverg(\VF)=p_t \VF_r^2 \varphi + \VF_r \varphi \diverg(p_t \VF_r)$, we deduce
that
\begin{equation}
  \int_{\MHyp} \VF_r^2 \varphi\,  p_t \dm  + \int_{\MHyp} \VF_r \varphi\, \diverg(p_t \VF_r) \dm
  = \int_{\dM} \VF_r \varphi\, p_t \scal{\VF_r}{\nu} \dsigma \,.
\end{equation}
This allows us to rewrite the left-hand side of \eqref{eq:ResultatIPP}
as
\begin{align}
  \int_{\M} \OpKolBack \varphi\, p_t\, \dm \;=\; & \int_{\MHyp} \VF_{0} \varphi\, p_t
  \dm \;+\; \frac{1}{2} \sum_{r=1}^{n} \int_{\MHyp} \VF_{r}^{2} \varphi\, p_t \dm
  \jbNoNum\\
  \;=\; & \int_{\MHyp} \left( \VF_{0} \varphi\, p_t - \frac{1}{2} \sum_{r=1}^{n}
    \VF_r \varphi\, \diverg(p_t \VF_r) \right) \dm
  \jbNoNum\\
  & + \frac{1}{2} \sum_{r=1}^{n} \int_{\dM} \VF_r \varphi\, p_t \scal{\VF_r}{\nu}
  \dsigma
  \jbNoNum\\
  \;=\; & \int_{\MHyp} J_t \varphi \dm + \frac{1}{2} \sum_{r=1}^{n} \int_{\dM} \VF_r
  \varphi\, p_t \scal{\VF_r}{\nu} \dsigma
  \label{eq:ApplicationDivergenceA}
\end{align}
We invoke the divergence formula~\eqref{eq:Stokes} once more, with
$A=\varphi J_t$, to deduce that
\begin{equation}
  \label{eq:ApplicationDivergenceB}
  \int_{\MHyp}  J_t \varphi \dm \;+\;  \int_{\MHyp} \varphi\, \diverg J_t\, \dm 
  \;=\;  \int_{\dM}\, \varphi\, \Jout_t\, \dsigma - \int_{H}\, \varphi\, \Jin_t\, \dsigma \,.
\end{equation}
Finally, combining equations
\eqref{eq:ApplicationDivergenceA}--\eqref{eq:ApplicationDivergenceB}
and using that $\diverg J_t = - \OpKolForw p_t$ yields the result.

\subsection{Proof of the extension lemma}
\label{appendix:ExtLemma}

Let $\mathscr{U}= \{ U_{\alpha},\, \alpha\in A \}$ be the set of all precompact
coordinate domains $U_{\alpha}$, with coordinate map $\Psi_{\alpha}: x\in U_{\alpha} \mapsto
(x^1, \ldots, x^n)\in\Rset^n$, such that:
\begin{enumerate}[a)\ ]
\item\label{Nb1} either $U_{\alpha} \cap \Hyp = \varnothing$ or $U_{\alpha} \cap \dM =
  \varnothing$;
\item\label{Nb2} if $U_{\alpha} \cap \dM \neq \varnothing$, there is $\varepsilon > 0$ such that
  \begin{align}
    \Psi_{\alpha} (U_{\alpha})        &\;=\; [0;\varepsilon)\times(-\varepsilon;\varepsilon)^{n-1} \,, \jbNoNum\\
    \Psi_{\alpha} (U_{\alpha} \cap \dM) &\;=\;  \{0\}\times(-\varepsilon;\varepsilon)^{n-1} \,,\jbNoNum
  \end{align}
  and $\frac{\partial}{\partial x^1}$ is the inward-pointing unit normal vector to
  $\dM$, i.e. $\frac{\partial}{\partial\nu} = - \frac{\partial}{\partial x^1}$;
\item\label{Nb3} if $U_{\alpha} \cap \Hyp \neq \varnothing$, there is $\varepsilon > 0$
  such that
  \begin{align}
    \Psi_{\alpha} (U_{\alpha})        &\;=\; (-\varepsilon;\varepsilon)^n \,, \jbNoNum\\
    \Psi_{\alpha} (U_{\alpha} \cap \dM) &\;=\;  \{0\}\times(-\varepsilon;\varepsilon)^{n-1} \,.\jbNoNum
  \end{align}
\end{enumerate}
Any $x\in\Mb$ has a precompact coordinate neighbourhood satisfying
either (\ref{Nb2}) or (\ref{Nb3}). It can be shrinked to satisfy
(\ref{Nb1}) since $H$ is a closed subset of $\M$ according to
assumption \ref{ass:Reset}.  Therefore, $\mathscr{U}$ is an open cover
of $\Mb$.

Let $\{ g_{\alpha},\, \alpha\in A \}$ be a partition of unity subordinate to
$\mathscr{U}$. For each $\alpha\in A$, we construct a function $\varphi_{\alpha}$ on
$U_{\alpha}$ as follows:
\begin{enumerate}[-]
\item If $\supp g_{\alpha}$ intersects neither $\dM$ nor $\Hyp$, we set
  $\varphi_{\alpha}=0$;
\item If $\supp g_{\alpha} \cap \dM \neq \varnothing$, we define $\varphi_{\alpha}$ in
  coordinates\footnote{We identify a point $x\in\Mb$ with its
    coordinates $(x^1, \ldots, x^n)$ in the chart $(U_{\alpha},\Psi_{\alpha})$.} by
  \begin{equation}
    \varphi_{\alpha} (x^1,x') \;=\; g_{\alpha} (x^1,x')\, 
    \left[
      \eta_1 (0,x') - x^1\, \eta_2 (0,x')
    \right] \,,
    \jbNoNum
  \end{equation}
  where $x'=(x^2, \ldots, x^n)$;
\item If  $\supp g_{\alpha} \cap \dM \neq \varnothing$, we set 
  \begin{equation}
    \varphi_{\alpha} (x^1,x') \;=\; g_{\alpha} (x^1,x')\; (\eta_1\circ\Phi^{-1}) (0,x') \,.
    \jbNoNum
  \end{equation}
\end{enumerate}
These definitions are compatible since, according to (\ref{Nb1}),
$U_{\alpha}$ cannot intersect both $\dM$ and $\Hyp$. $\varphi_{\alpha}$ is of class
$\Cclass^2$ since the $\eta_i$'s are $\Cclass^2$ by hypothesis and $\phi$ is
a $\Cclass^2$-diffeomorphism by \ref{ass:Reset}.  Moreover, $\varphi_{\alpha}$ is
compactly supported in $U_{\alpha}$ because $\supp \varphi_{\alpha} \subset \supp g_{\alpha} \subset
U_{\alpha}$. We extend $\varphi_{\alpha}$ to $\Stb$ by $\varphi_{\alpha}=0$ outside $U_{\alpha}$. We
observe that $\varphi_{\alpha}\neq0$ if and only if $\supp g_{\alpha}$ intersects either
$\supp \eta_1 \cup\, \supp \eta_2$ or $\Phi(\supp \eta_1 \cap G)$, which are both
compact sets. Consequently, $\{\, \alpha \,|\, \varphi_{\alpha} \neq 0 \,\} $ is a finite
set because $\{ \supp g_{\alpha}, \alpha \in A \}$ is locally finite.

Finally, we set $\varphi=\sum_{\alpha} \varphi_{\alpha}$. This is in fact a finite sum by the
preceding observation, which implies that $\phi \in \Ckc{2}{\Stb}$.  The
properties (\ref{it:L1})-(\ref{it:L3}) of the lemma follow then
easily from our construction, using that $\sum_{\alpha} g_{\alpha}=1$, and the
proof of Lemma~\ref{lem:ExtensionLemma} is complete.

\bibliographystyle{elsart-num}
\bibliography{fpke}

\end{document}